\documentclass[12pt]{amsart}
\usepackage{amssymb,amsmath, amsthm,latexsym}
\usepackage{a4wide}

\theoremstyle{plain}

\newtheorem{lem}{Lemma}[section]
\newtheorem{theo}[lem]{Theorem}
\newtheorem{prop}[lem]{Proposition}
\newtheorem{corollary}[lem]{Corollary}

\parskip=\medskipamount
\font\k=cmr7
\font\rm=cmr12

  \newcommand {\cu}{\mbox{\k cus}}
  \newcommand {\di}{\mbox{\k dis}}
  \newcommand {\ac}{\mbox{\k c}}
  
  \newcommand {\res}{\mbox{\k res}}

  \newcommand {\temp}{\mbox{\k temp}}
  \newcommand {\comp}{\mbox{\k comp}}
  \newcommand {\C}{{\mathbb C}}
  \newcommand {\bH}{{\mathbb H}}
  \newcommand {\N}{{\mathbb N}}
  \newcommand {\R}{{\mathbb R}}
  \newcommand {\Z}{{\mathbb Z}}
  \newcommand {\Q}{{\mathbb Q}}

  \newcommand {\af}{{\mathfrak a}}

\renewcommand {\H}{{\mathcal H}}

  \newcommand {\G}{{\bf G}}
  \newcommand {\bP}{{\bf P}}

  \newcommand {\E}{{\mathcal E}}
  \newcommand {\cD}{{\mathcal D}}

 \newcommand {\cA}{{\mathcal A}}
\newcommand {\cT}{{\mathcal T}}

\newcommand {\bs}{\backslash}
\newcommand {\ds}{/\hskip-2pt/}

\renewcommand{\Im}{\operatorname{Im}}
\renewcommand{\Re}{\operatorname{Re}}
\newcommand{\Tr}{\operatorname{Tr}}
\newcommand{\End}{\operatorname{End}}
\newcommand{\sign}{\operatorname{sign}}

\newcommand{\Id}{\operatorname{Id}}

\newcommand{\rk}{\operatorname{rank}}
\newcommand{\vol}{\operatorname{vol}}
\newcommand{\Area}{\operatorname{Area}}

\newcommand{\SL}{\operatorname{SL}}
\newcommand{\GL}{\operatorname{GL}}
\newcommand{\SO}{\operatorname{SO}}
\newcommand{\PSL}{\operatorname{PSL}}

\newcommand{\supp}{\operatorname{supp}}
\newcommand{\sub}{\operatorname{sub}}

\newcommand{\diver}{\operatorname{div}}
\newcommand{\grad}{\operatorname{grad}}

\begin{document}
\title[Weyl law]
{Weyl's law in the theory of automorphic forms}
\date{\today}

\author{Werner M\"uller}
\address{Universit\"at Bonn\\
Mathematisches Institut\\
Beringstrasse 1\\
D -- 53115 Bonn, Germany}
\email{mueller@math.uni-bonn.de}
\keywords{spectrum, automorphic forms}
\subjclass{Primary: 58G25, Secondary: 22E40}
\begin{abstract}
For a compact Riemannian manifold, Weyl's law describes the asymptotic 
behavior of the counting function of the eigenvalues of the associated Laplace
operator.
In this paper we discuss Weyl's law in the context of automorphic forms. 
The underlying manifolds are locally symmetric spaces of finite volume. In the
non-compact case Weyl's law is closely related to the problem of
existence of cusp forms. 
\end{abstract}
\maketitle

\section{Introduction}
Let $M$ be a smooth, compact Riemannian manifold  of dimension $n$ with 
smooth boundary $\partial M$ (which may be empty). Let 
\[\Delta=-\diver\circ\grad=d^\ast d\]
be the Laplace-Beltrami operator associated with the metric $g$ of $M$.
We consider the Dirichlet eigenvalue problem
\begin{equation}\label{1.1}
\Delta\phi=\lambda\phi,\quad \phi\big|_{\partial M}=0.
\end{equation}
As is well known, (\ref{1.1}) has a discrete set of solutions
\[0\le\lambda_0\le\lambda_1\le\cdots\to\infty\]
whose only accumulation point is at infinity and each eigenvalue occurs
with finite multiplicity. The corresponding eigenfunctions $\phi_i$
can be chosen such that $\{\phi_i\}_{i\in\N_0}$ is an orthonormal basis
of $L^2(M)$. A fundamental problem in analysis on manifolds is to study  
the distribution of the eigenvalues of $\Delta$ and their relation to 
the geometric
and topological structure of the underlying manifold. One of the first results
in this context is Weyl's law for the asymptotic behavior of the eigenvalue
counting function. For $\lambda\ge0$ let
\[N(\lambda)=\#\bigl\{j\colon \sqrt{\lambda_j}\le\lambda\big\}\]
be the counting function of the eigenvalues of $\sqrt{\Delta}$,
where eigenvalues are counted with  multiplicities. Denote by 
${\bf \Gamma}(s)$ the Gamma function.
Then the Weyl law states
\begin{equation}\label{1.2}
N(\lambda)=\frac{\vol(M)}{(4\pi)^{n/2}{\bf \Gamma}\left(\frac{n}{2}+1\right)}
\lambda^{n}+o(\lambda^{n}),\quad \lambda\to\infty.
\end{equation}
This was first proved by Weyl \cite{We1} for a bounded domain $\Omega\subset
\R^3$. Written in a slightly different form it is known in
physics as the Rayleigh-Jeans law. Raleigh \cite{Ra} derived it for a cube.  
Garding \cite{Ga}  proved Weyl's law for a general elliptic operator on a 
domain in $\R^n$.
For a closed Riemannian manifold (\ref{1.2}) was proved by Minakshisundaram
and Pleijel \cite{MP}. 

Formula (\ref{1.2}) does not say very much about the finer structure of the
eigenvalue distribution. The basic question is
 the estimation of the remainder term
\[R(\lambda):=N(\lambda)
-\frac{\vol(M)}{(4\pi)^{n/2}{\bf \Gamma}\left(\frac{n}{2}+1\right)}
\lambda^{n}.\] 
That this is a deep problem shows the following example. Consider the flat
2-dimensional torus $T=\R^2/(2\pi \Z)^2$. Then the eigenvalues of the flat 
Laplacian
are $\lambda_{m,n}:=m^2+n^2$, $m,n\in\Z$ and the counting function
equals
\[N(\lambda)=\#\big\{(m,n)\in\Z^2\colon \sqrt{m^2+n^2}\le \lambda\big\}.\]
Thus $N(\lambda)$ is the number of lattice points in the circle of radius
$\lambda$. An elementary packing argument, attributed to Gauss, gives
\[N(\lambda)=\pi\lambda^2+O\left(\lambda\right),\]
and the circle problem is to find the best exponent $\mu$ such that
\[N(\lambda)=\pi\lambda^2+O_\varepsilon\left(\lambda^{\mu+\varepsilon}\right),
\quad \forall \varepsilon>0.\]
The conjecture of Hardy is $\mu=1/2$. The first nontrivial result is due to 
Sierpinski who showed that one can take $\mu=2/3$. Currently the best known 
result is $\mu=131/208\approx 0.629$ which is due to Huxley. 
Levitan \cite{Le} has shown that for a domain in $\R^n$ the remainder term is
of order $O(\lambda^{n-1})$.
 
For a closed Riemannian manifold, Avakumovi\'c \cite{Av} proved the Weyl
estimate with optimal remainder term:
\begin{equation}\label{1.3}
N(\lambda)=\frac{\vol(M)}{(4\pi)^{n/2}{\bf \Gamma}\left(\frac{n}{2}+1\right)}
\lambda^{n}+O(\lambda^{n-1}),\quad \lambda\to\infty.
\end{equation}
This result was extended to more general, and higher order operators by
H\"ormander \cite{Ho}. 
As shown by Avakumovi\'c the bound $O(\lambda^{n-1})$ of the remainder term
is optimal for the sphere. On the other hand, under certain assumption on the 
geodesic flow, the estimate can be slightly improved. Let $S^\ast M$ be the
unit cotangent bundle and let $\Phi_t$ be the geodesic flow. Suppose that
the set of $(x,\xi)\in S^\ast M$ such that $\Phi_t$ has a contact of infinite
order with the identity at $(x,\xi)$ for some $t\not=0$, has measure zero in
$S^\ast M$. Then Duistermaat and Guillemin \cite{DG} 
proved that the remainder term satisfies $R(\lambda)=o(\lambda^{n-1})$. This is
a slight improvement over (\ref{1.3}).

In \cite{We3} Weyl formulated a  conjecture which claims the existence of a
second term in the asymptotic expansion for a bounded domain $\Omega\subset
\R^3$, namely he predicted that
\[ N(\lambda)=\frac{\vol(\Omega)}{6\pi^2}\lambda^{3}-\frac{\vol(\partial
\Omega)}{16\pi}\lambda^2+o(\lambda^2).\]
This was proved for manifolds with boundary under a certain condition on the 
periodic billiard trajectories, by Ivrii \cite{Iv} and Melrose \cite{Me}.

The purpose of  this paper is to discuss Weyl's law in the context of locally 
symmetric spaces $\Gamma\bs S$ of finite volume and non-compact type.
Here $S=G/K$ is a Riemannian
symmetric space, where $G$ is  a real semi-simple Lie group  
of non-compact type, and $K$ a maximal compact subgroup of $G$. Moreover
$\Gamma$ is a lattice in $G$, i.e., a discrete subgroup of finite 
covolume. Of particular interest are 
arithmetic subgroups such as  the principal
congruence subgroup  $\Gamma(N)$ of $\SL(2,\Z)$ of level $N\in\N$. 
Spectral theory of the Laplacian on arithmetic quotients
$\Gamma\bs S$ is intimately related with the theory of automorphic forms.
In fact, for a symmetric space $S$ it is more natural and important 
to consider not only the Laplacian,
but the whole algebra $\cD(S)$ of $G$-invariant differential operators on 
$S$. It is known that $\cD(S)$ is a finitely generated commutative algebra 
\cite{He}. Therefore, it makes sense to study the joint spectral decomposition 
of $\cD(S)$. 
Square integrable joint eigenfunctions of $\cD(S)$  are examples of 
automorphic forms. Among them
are the cusp forms which satisfy additional decay conditions. Cusps forms are
 the building blocks of the theory of automorphic forms and, according to deep
and far-reaching conjectures of Langlands \cite{La2}, are expected to
provide important relations between harmonic analysis  and number theory.

Let $G=NAK$ be the Iwasawa decomposition of $G$ and let $\af$ be the Lie 
algebra of $A$.   
If $\Gamma\bs S$ is compact, the spectrum of $\cD(S)$ in $L^2(\Gamma\bs S)$
is a discrete subset of the complexification $\af^\ast_\C$ of $\af^\ast$. It
 has been studied by Duistermaat, Kolk, and Varadarajan in 
\cite{DKV}. The method is based on the Selberg trace formula. The results
are more refined statements about the distribution of the spectrum than
just the Weyl law. For example, one gets  estimations for the distribution of
the tempered and the complementary spectrum. We will  review briefly these 
results in section 2.  

If $\Gamma\bs S$ is non-compact, which is the case for many important 
arithmetic groups,  
the Laplacian has a large continuous spectrum which can be 
described in terms of Eisenstein series \cite{La1}. Therefore, it is 
not obvious that the Laplacian has any eigenvalue $\lambda>0$, and an important
problem in the theory of automorphic forms is the existence and construction
of cusp forms for a given lattice $\Gamma$. This is were the Weyl law comes 
into play. Let $\bH$ be the upper half-plane. Recall that $\SL(2,\R)$ acts on
$\bH$ by fractional linear transformations.
Using his trace formula \cite{Se2}, Selberg established the 
following version of Weyl's law for an arbitrary lattice $\Gamma$ in 
$\SL(2,\R)$
\begin{equation}\label{1.4}
N_\Gamma(\lambda)+M_\Gamma(\lambda)\sim\frac{\Area(\Gamma\bs\bH)}
{4\pi}\lambda^2,\quad \lambda\to\infty
\end{equation}
\cite[p. 668]{Se2}.
Here $N_\Gamma(\lambda)$ is the counting function of the eigenvalues and
 $M_\Gamma(\lambda)$ is the winding number of the determinant 
$\phi(1/2+ir)$ of the scattering matrix which is given by the
constant Fourier coefficients of the Eisenstein series (see section 4). In
general, the two functions on the left can not be estimated separately. 
However, for congruence groups like $\Gamma(N)$, the meromorphic function 
$\phi(s)$ can be expressed in 
terms of well-known functions of analytic number theory. In this case, it is
possible to show that the growth of $M_\Gamma(\lambda)$ is of lower order 
which implies Weyl's law for the counting function of the eigenvalues
\cite[p.668]{Se2}.
Especially it follows  that Maass
cusp forms exist in abundance for congruence groups. 
On the other hand, there are
indications \cite{PS1}, \cite{PS2} that the existence of many  cusp forms
may be restricted to arithmetic groups. This will be discussed in detail
in section 4. 

In section 5 we discuss the general case of a non-compact arithmetic quotient 
$\Gamma\bs S$. There has been some recent progress with the spectral
 problems discussed above. 
Lindenstrauss and Venkatesh \cite{LV} established Weyl's law 
without  remainder term   for congruence subgroups of a
split adjoint semi-simple group $\G$. In \cite{Mu3} this had been proved for 
congruence subgroups of $\SL(n)$ and for the Bochner-Laplace operator acting
in sections of a locally homogeneous vector bundle over $S_n=\SL(n,\R)/\SO(n)$.
For  congruence subgroups of $\G=\SL(n)$,  an estimation of 
the remainder term in Weyl's law  has been established by E. Lapid and the 
author in \cite{LM}. 
Using the approach of \cite{DKV}
combined with the Arthur trace formula,
the results of \cite{DKV} have been extended  in \cite{LM} to the cuspidal 
spectrum of $\cD(S_n)$.

\section{Compact locally symmetric spaces}
\setcounter{equation}{0} 
In this section we review H\"ormanders method of the derivation of Weyl's law
with remainder term for the Laplacian $\Delta$ of a closed Riemannian manifold
$M$ of dimension $n$. Then we will discuss the results of \cite{DKV}
concerning spectral asymptotics for  compact locally symmetric manifolds.

The method of H\"ormander \cite{Ho}  to estimate the remainder term
is based on the study of the kernel of $e^{-it\sqrt{\Delta}}$. 
The main point is the construction of a good approximate fundamental solution
to the wave equation by means of the theory of Fourier integral operators and
the analysis of the singularities of its trace
\[\Tr e^{-it\sqrt{\Delta}}=\sum_j e^{-it\sqrt{\lambda_j}},\]
which is well-defined as a distribution.  The analysis of H\"ormander of 
the ``big'' singularity of  
$\Tr e^{-it\sqrt{\Delta}}$ at $t=0$ leads to the following key 
result \cite[(2.16)]{DG}. Let $\mu_j:=\sqrt{\lambda_j}$, $j\in\N$.
There exist $c_j\in\R$,
$j=0,...,n-1$, and $\varepsilon>0$ such that for every 
$h\in\mathcal{S}(\R)$ with 
$\supp\hat h\subset[-\varepsilon,\varepsilon]$ and $\hat h\equiv 1$ in
a neighborhood of $0$ one has
\begin{equation}\label{2.1}
\sum_j h(\mu-\mu_j)\sim(2\pi)^{-n}
\sum_{k=0}^{n-1} c_k\mu^{n-1-k},\quad \mu\to\infty,
\end{equation}
and rapidly decreasing as $\mu\to-\infty$. The constants $c_k$ are of the form
\[c_k=\int_M\omega_k,\]
where the $\omega_k$'s are real valued smooth densities on $M$ canonically 
associated to the Riemannian metric of $M$. Especially
\[c_0=\vol(S^\ast M),\quad c_1=(1-n)\int_{S^\ast M}\sub \Delta,\]
where $S^\ast M$ is the unit co-tangent bundle, and $\sub\Delta$ denotes the
subprincipal symbol of $\Delta$. Consideration of the top term in (\ref{2.1})
leads to the basic estimates for the eigenvalues.

If $M=\Gamma\bs G/K$ is a locally symmetric manifold, the Selberg trace formula
can be used to replace (\ref{2.1}) by an exact formula \cite{DKV}. Actually,
if the rank of $M$ is bigger than 1, the spectrum is multidimensional. Then the
Selberg trace formula gives more refined information.

As example,  we consider a compact hyperbolic surface 
$M=\Gamma\bs\bH$, where $\Gamma\subset\PSL(2,\R)$ is a discrete, torsion-free,
co-compact subgroup. Let $\Delta$ be the hyperbolic Laplace operator which
is given by
\begin{equation}\label{2.1b}
\Delta=-y^2\left(\frac{\partial^2}{\partial x^2}+\frac{\partial^2}{\partial y^2}
\right),\quad z=x+iy.
\end{equation}
Write 
the eigenvalues $\lambda_j$ of $\Delta$ as 
\[\lambda_j=\frac{1}{4}+r_j^2,\]
where $r_j\in\C$ and $\arg(r_j)\in\{0,\pi/2\}$.
Let $h$ be an analytic function in a strip
$|\Im(z)|\le\frac{1}{2}+\delta$, $\delta>0$, such that
\begin{equation}\label{2.1a}
h(z)=h(-z),\quad |h(z)|\le C(1+|z|)^{-2-\delta}.
\end{equation}
Let 
$$g(u)=\frac{1}{2\pi}\int_\R h(r)e^{iru}dr.$$
Given $\gamma\in\Gamma$ denote by $\{\gamma\}_\Gamma$ its 
$\Gamma$-conjugacy class.
Since $\Gamma$ is co-compact, each $\gamma\not= e$ is hyperbolic. Each 
hyperbolic element 
$\gamma$ is the power of a primitive hyperbolic element $\gamma_0$.
A hyperbolic
conjugacy class determines a closed geodesic $\tau_\gamma$ of $\Gamma\bs\bH$. 
Let $l(\gamma)$ denote the length of $\tau_\gamma$. 
Then the Selberg trace formula \cite{Se1} is the following identity:
\begin{equation}\label{2.2}
\sum_{j=0}^\infty h(r_j)=\frac{\Area(\Gamma\bs\bH)}{4\pi}\int_\R h(r)r
\tanh(\pi r)\;dr
+\sum_{\{\gamma\}_\Gamma\not=e} 
\frac{l(\gamma_0)}{2\sinh \left( \frac{l(\gamma)}{2}\right)} g(l(\gamma)).
\end{equation}
Now let $g\in C^\infty_c(\R)$ and $h(z)=\int_\R g(r)e^{-ir z}\;dr$. Then $h$ is
entire and rapidly decreasing in each strip $|\Im(z)|\le c$, $c>0$. Let
$t\in\R$ and set
\[h_t(z)=h(t-z)+h(t+z).\]
Then $h_t$ is entire and satisfies (\ref{2.1a}). Note that 
$\hat h_t(r)=e^{-itr}g(r)+e^{itr}g(-r).$
We symmetrize the spectrum by $r_{-j}:=-r_j$, $j\in\N$. 
Then by (\ref{2.2}) we get
\begin{equation}\label{2.3}
\begin{split}
\sum_{j=-\infty}^\infty h(t-r_j)=&\frac{\Area(\Gamma\bs\bH)}{2\pi}\int_\R h(t-r)r
\tanh(\pi r)\;dr\\
&+\sum_{\{\gamma\}_\Gamma\not=e} 
\frac{l(\gamma_0)}{2\sinh \left( \frac{l(\gamma)}{2}\right)} \left(
e^{-itl(\gamma)}g(l(\gamma))+e^{itl(\gamma)}g(-l(\gamma))\right).
\end{split}
\end{equation}
Let $\varepsilon>0$ be such that $l(\gamma)>\varepsilon$ for all hyperbolic
conjugacy classes $\{\gamma\}_\Gamma$.  The following lemma is an immediate 
consequence of (\ref{2.3}).
\begin{lem}\label{l2.1} Let $g\in C_c^\infty(\R)$ such that 
$\supp g\subset (-\varepsilon,\varepsilon)$. 
Let $h(z)=\int_\R g(r)e^{-ir z}\;dr$. Then for all $t\in\R$ we
have
\begin{equation}\label{2.4a}
\sum_{j=-\infty}^\infty h(t-r_j)=\frac{\Area(\Gamma\bs\bH)}{2\pi}
\int_\R h(t-r)r\tanh(\pi r)\;dr.
\end{equation}
\end{lem} 
Changing variables in the integral on the right and using that
\[\tanh(\pi(r+t))=1-\frac{2e^{-2\pi(r+t)}}{1+e^{-2\pi(r+t)}}=
-1+\frac{2e^{2\pi(r+t)}}{1+e^{2\pi(r+t)}},\]
we obtain the following asymptotic expansion
\begin{equation}\label{2.5}
\begin{split}
\sum_{j=-\infty}^\infty h(t-r_j)=
\frac{\Area(\Gamma\bs\bH)}{2\pi}\left(|t|\int_\R h(r)\;dr-
\sign t\int_\R h(r)r\;dr\right)
+O\left(e^{-2\pi|t|}\right),
\end{split}
\end{equation}
as $|t|\to\infty$. If $h$ is even, the second term vanishes and the 
asymptotic expansion is related to (\ref{2.1}).
The asymptotic expansion (\ref{2.5})
can be used to derive  estimates for the number of eigenvalues near a given 
point $\mu\in\R$. 
\begin{lem}\label{l2.2} 
For every $a>0$ there exists $C>0$ such that 
$$\#\{j\colon |r_j-\mu|\le a\}\le C(1+|\mu|)$$
for all $\mu\in\R$.
\end{lem}
\begin{proof}  We proceed as in the proof of Lemma 2.3 in \cite{DG}. 
 As shown in the proof, there exists $h\in\mathcal{S}(\R)$
such that $h\ge 0$, $h>0$ on $[-a,a]$, $\hat h(0)=1$,
and $\supp \hat h$ is contained in any prescribed neighborhood of $0$. 
Now observe that there are only finitely many eigenvalues $\lambda_j=1/4+r_j^2$
with $r_j\notin\R$. Therefore it suffices to consider $r_j\in\R$. Let 
$\mu\in\R$. By  (\ref{2.5})  we get
$$\#\{j\colon |r_j-\mu|\le a,\;r_j\in\R\}\cdot 
\min\{h(u)\colon |u|\le a\}\le
\sum_{r_j\in\R} h(\mu-r_j)\le C(1+|\mu|).$$
\end{proof}
This lemma is the basis of the following auxiliary results.
\begin{lem}\label{l2.3}
For every $h$ as above there exists $C>0$ such that
\begin{equation}\label{1.9}
\sum_{|r_j|\le\lambda}\bigg|\int_{\R-[-\lambda,\lambda]} h(t-r_j)\;dt\bigg|\le
C\lambda,\quad
\sum_{|r_j|>\lambda}\bigg|\int_{-\lambda}^\lambda h(t-r_j)\;dt\bigg|\le
C\lambda,
\end{equation}
for all $\lambda\ge 1$.
\end{lem}
\begin{proof}
Since $h$ is rapidly decreasing, there exists $C>0$ such
 that 
$|h(t)|\le C (1+|t|)^{-4}$, $t\in\R$.
Let $[\lambda]$ be the largest integer $\le\lambda$. Then we get
\begin{equation*}
\begin{split}
\sum_{|r_j|\le\lambda}\bigg|\int_\lambda^\infty &h(t-r_j)\;dt\bigg|\le
\sum_{|r_j|\le\lambda}\int_{\lambda-r_j}^\infty |h(t)|\;dt
\le C \sum_{|r_j|\le\lambda}\frac{1}{(1+\lambda-r_j)^3}\\
&=\sum_{k=-[\lambda]}^{[\lambda]-1}\sum_{k\le r_j\le k+1}
\frac{1}{(1+\lambda-r_j)^3}\le \sum_{k=-[\lambda]}^{[\lambda]-1}
\frac{\#\{j\colon |r_j-k|\le 1\}}{(\lambda-k)^3},
\end{split}
\end{equation*}
and by Lemma \ref{l2.2} the right hand side is  bounded by $C\lambda$ for
$\lambda\ge 1$.
Similarly we get
$$\sum_{|r_j|\le\lambda}\bigg|\int^{-\lambda}_{-\infty} h(t-r_j)\;dt\bigg|\le
C_2\lambda.$$
The second series can be treated in the same way.
\end{proof}
\begin{lem}\label{l2.4}
Let $h$ be as in Lemma \ref{l2.1} and such that $\hat h(0)=1$. Then 
\begin{equation}\label{1.10}
\int_{-\lambda}^\lambda \sum_{j=-\infty}^\infty h(t-r_j)\;dt=
\frac{\Area(\Gamma\bs\bH)}{2\pi}\lambda^2+O(\lambda)
\end{equation}
as $\lambda\to\infty$.
\end{lem}
\begin{proof}
To prove the lemma, we integrate (\ref{2.4a}) and  determine
 the asymptotic behavior of the integral on the right. Let
 $p(r)$ be a continuous function on $\R$ such that $|p(r)|\le C(1+|r|)$
and $p(r)=p(-r)$.
Changing the order of integration and using that $\int_\R h(t-r)\;dt=\hat 
h(0)=1$,  we get
\begin{equation*}
\begin{split}
\int_{-\lambda}^\lambda \int_\R h(t-r)p(r)\;dr\;dt&=\int_{-\lambda}^\lambda p(r)\;dr
-\int_{-\lambda}^\lambda \left(\int_{\R-[-\lambda,\lambda]}
h(t-r)\;dt\right)p(r)\;dr\\
&\quad + \int_{\R-[-\lambda,\lambda]}\left(\int_{-\lambda}^\lambda  h(t-r)\;dt
\right)p(r)\;dr.
\end{split}
\end{equation*}
Let $C_1>0$ be such that $|h(r)|\le C_1(1+|r|)^{-3}$, $r\in\R$. Then the the
second and the third integral can be estimated by $C(1+\lambda)$. 
Thus we get
\begin{equation}\label{2.9}
\int_{-\lambda}^\lambda \int_\R h(t-r)p(r)\;dr\;dt=\int_{-\lambda}^\lambda
p(r)\;dr+O(\lambda),\quad \lambda\to\infty.
\end{equation}
If we apply (\ref{2.9}) to  $p(r)=r\tanh(\pi r)$, we obtain
\begin{equation}\label{2.10}
\int_{-\lambda}^\lambda \int_\R h(t-r)r\tanh(\pi r)\;dr\;dt= \lambda^2+
O(\lambda).
\end{equation}
This proves the lemma.
\end{proof}
We are now ready to prove Weyl's law.
We choose $h$ such that $\hat h$ has sufficiently small  support and 
$\hat h(0)=1$. Then
\begin{equation*}
\begin{split}
\int_{-\lambda}^\lambda\sum_{j=-\infty}^\infty h(t-r_j)\;dt
&=\sum_{|r_j|\le\lambda}\int_\R h(t-r_j)\;dt-\sum_{|r_j|\le\lambda}
\int_{\R-[-\lambda,\lambda]} h(t-r_j)\;dt\\
&\quad +\sum_{|r_j|>\lambda}
\int_{-\lambda}^\lambda h(t-r_j)\;dt.
\end{split}
\end{equation*}
Using that $\int_\R h(t-r)\;dt=\hat h(0)=1$, we get
\begin{equation*}
\begin{split}
2 N_\Gamma(\lambda)=\int_{-\lambda}^\lambda\sum_j h(t-r_j)\;dt
&+\sum_{|r_j|\le\lambda}
\int_{\R-[-\lambda,\lambda]} h(t-r_j)\;dt\\
&\quad-\sum_{|r_j|>\lambda}
\int_{-\lambda}^\lambda h(t-r_j)\;dt.
\end{split}
\end{equation*}
By Lemmas \ref{l2.3} and  \ref{l2.4} we obtain
\begin{equation}\label{2.10a}
N_\Gamma(\lambda)=\frac{\Area(\Gamma\bs\bH)}{4\pi}\lambda^2+O(\lambda).
\end{equation}

We turn now to an arbitrary Riemannian symmetric space $S=G/K$ of non-compact
type and we review the main results of \cite{DKV}.
The group of motions  $G$ of $S$ is a semi-simple Lie group of non-compact 
type with finite center and $K$ is a maximal compact subgroup of $G$. 
 The Laplacian $\Delta$ of $S$ is a $G$-invariant differential operator on 
$S$, i.e., $\Delta$ commutes with the left translations $L_g$, $g\in G$. 
Besides of $\Delta$ we need to consider the ring $\cD(S)$ 
of all invariant 
differential operators on $S$. It is well-known that  $\cD(S)$ is commutative
and finitely generated. Its structure can be described as follows. 
Let $G=NAK$ be the Iwasawa decomposition of $G$, $W$ the Weyl group of $(G,A)$
 and  $\af$ be the Lie algebra of $A$.  
Let $S(\af_\C)$ be the symmetric algebra of the complexification
$\af_\C=\af\otimes\C$ of $\af$ and let $S(\af_\C)^W$ be the subspace of
Weyl group invariants in $S(\af_\C)$. Then by a theorem of
Harish-Chandra \cite[Ch. X, Theorem 6.15]{He} there is a canonical isomorphism
\begin{equation}\label{2.11}
\mu\colon\cD(S)\cong S(\af_\C)^W.
\end{equation}
This shows that $\cD(S)$ is commutative. The minimal number of generators
equals the rank of $S$ which is $\dim\af$ \cite[Ch.X, \S 6.3]{He}. Let 
$\lambda\in\af_\C^\ast$. Then by (\ref{2.11}), $\lambda$ determines an 
 character
\[\chi_\lambda\colon\cD(S)  \to \C\]
and  $\chi_\lambda=\chi_{\lambda^\prime}$ if and only if $\lambda$ and 
$\lambda^\prime$ are in the same $W$-orbit.
Since  $S(\af_\C)$ is integral over $S(\af_\C)^W$ \cite[Ch. X, Lemma 6.9]{He},
each  character of $\cD(S)$ is of the form $\chi_\lambda$ for 
some $\lambda\in\af_\C^\ast$. Thus the characters of $\cD(S)$ are
parametrized by $\af_\C^\ast/W$.

Let $\Gamma\subset G$ be a discrete, torsion-free, co-compact subgroup of $G$. 
Then $\Gamma$ acts properly discontinuously on $S$ without fixed points
and the quotient $M=\Gamma\bs S$ is a locally symmetric manifold which is
equipped with the  metric induced from the invariant metric of $S$. Then each 
$D\in\cD(S)$ descends to a differential operator 
\[D\colon C^\infty(\Gamma\bs S)\to C^\infty(\Gamma\bs S).\]
Let $\E\subset C^\infty(\Gamma\bs S)$ be an eigenspace of the Laplace
operator. Then $\E$ is a finite-dimensional vector
space which is invariant under $D\in\cD(S)$. For each $D\in\cD(S)$, the formal 
adjoint $D^\ast$ of $D$ also belongs to $\cD(S)$. Thus we get a representation
\[\rho\colon \cD(S)\to \End(\E)\]
by commuting normal operators. Therefore, $\E$ decomposes into the direct
sum of joint eigenspaces of $\cD(S)$. Given $\lambda\in\af_\C^\ast/W$, let 
\[\E(\lambda)=\{\varphi\in C^\infty(\Gamma\bs S)\colon D\varphi=\chi_\lambda(D)
\varphi,\;D\in\cD(S)\}.\]
Let $m(\lambda)=\dim\E(\lambda)$. Then the spectrum $\Lambda(\Gamma)$
of $\Gamma\bs S$ is defined to be
\[\Lambda(\Gamma)=\{\lambda\in\af_\C^\ast/W\colon m(\lambda)>0\},\]
and we get an orthogonal direct sum decomposition
\[L^2(\Gamma\bs S)=\bigoplus_{\lambda\in\Lambda(\Gamma)}\E(\lambda).\]

If we pick a fundamental domain for $W$,  we may regard $\Lambda(\Gamma)$ as a 
subset of $\af_\C^\ast$. If $\rk(S)>1$, then $\Lambda(\Gamma)$ is 
multidimensional. Again the distribution of $\Lambda(\Gamma)$ is studied using
the Selberg trace formula \cite{Se1}. To describe it we need to introduce some
notation. Let $C_c^\infty(G\ds K)$ be the subspace of all $f\in C^\infty_c(G)$
which are $K$-bi-invariant.  Let
\[\cA\colon C_c^\infty(G\ds K)\to C_c^\infty(A)^W\]
be the Abel transform which is defined by
\[\cA(f)(a)=\delta(a)^{1/2}\int_N f(an)\,dn,\quad a\in A,\]
where $\delta$ is the modulus function of the minimal parabolic subgroup
$P=NA$. Given $h\in C_c^\infty(A)^W$, let
\[\hat h(\lambda)=\int_A h(a)e^{\langle\lambda,H(a)\rangle}\;da.\]
Let $\beta(i\lambda)$, $\lambda\in\af^\ast$, be the Plancherel density.
Then the Selberg trace formula is the following identity
\begin{equation}\label{2.16}
\begin{split}
\sum_{\lambda\in\Lambda(\Gamma)}m(\lambda)\hat h(\lambda)=&
\frac{\vol(\Gamma\bs G)}{|W|}\int_{\af^\ast}\hat h(\lambda)\beta(i\lambda)
\;d\lambda\\
&+\sum_{[\gamma]_\Gamma\not=e}\vol(\Gamma_\gamma\bs G_\gamma)\int_{G_\gamma\bs G}
\cA^{-1}(h)(x^{-1}\gamma x)\;d_\gamma\bar x.
\end{split}
\end{equation} 
This  is still not the final form of the Selberg trace formula. The 
distributions 
\begin{equation}\label{2.16a}
J_\gamma(f)=\vol(\Gamma_\gamma\bs G_\gamma)\int_{G_\gamma\bs G}
f(x^{-1}\gamma x)\;d_\gamma\bar x, \quad f\in C^\infty_c(G),
\end{equation}
are invariant distribution an $G$ 
and can be computed using Harish-Chandra's Fourier inversion formula. This 
brings (\ref{2.16}) into a form which is similar to (\ref{2.2}). 
For the present purpose, however, it suffices to work with (\ref{2.16}). 
Since for $\gamma\not=e$,
the conjugacy class  of $\gamma$ in $G$ is closed and does not
intersect $K$,  there exists an open neighborhood $V$ of $1$ in 
$A$ satisfying $V=V^{-1}$, $V$ is invariant under $W$, and $J_\gamma(\cA^{-1}(h))
=0$ for all $h\in C^\infty_c(V)$ \cite[Propostion 3.8]{DKV}. Thus we get
\begin{equation}\label{2.17}
\sum_{\lambda\in\Lambda(\Gamma)}m(\lambda)\hat h(\lambda)=
\frac{\vol(\Gamma\bs G)}{|W|}\int_{\af^\ast}\hat h(\lambda)\beta(i\lambda)
\;d\lambda
\end{equation} 
for all $h\in C_c^\infty(V)$. One can now proceed as in the case of the upper
half-plane. The basic step is again to estimate the number of $\lambda\in
\Lambda(\Gamma)$ lying in a ball of radius $r$ around a variable point 
$\mu\in i\af^\ast$. This can be achieved by inserting appropriate test functions
$h$ into (\ref{2.17}) \cite[section 7]{DKV}. 
Let 
\[\Lambda_{\temp}(\Gamma)=\Lambda(\Gamma)\cap i\af^\ast,\quad 
\Lambda_{\comp}(\Gamma)=\Lambda(\Gamma)\setminus \Lambda_{\temp}(\Gamma)\]
be the tempered and complementary spectrum, respectively. Given an open
bounded subset $\Omega$ of $\af^\ast$ and $t>0$,  let 
\begin{equation}\label{2.17a}
t\Omega:=\{t\mu\colon \mu\in\Omega\}.
\end{equation}
One of the main results of \cite{DKV} is the following asymptotic formula
for the distribution of the tempered spectrum \cite[Theorem 8.8]{DKV}
\begin{equation}\label{2.18}
\sum_{\lambda\in\Lambda_{\temp}(\Gamma)\cap(it\Omega)}m(\lambda)=
\frac{\vol(\Gamma\bs G)}{|W|}\int_{it\Omega}\beta(i\lambda)\,d\lambda
+O(t^{n-1}),\quad t\to\infty,
\end{equation}
Note that the leading term is of order $O(t^n)$.
The growth of the complementary spectrum is of lower order.
Let $B_t(0)\subset \af^\ast_\C$ be the ball of radius $t>0$ around 0. 
There exists $C>0$ such that
for all $t\ge 1$  
\begin{equation}\label{2.19}
\sum_{\lambda\in\Lambda_{\comp}(\Gamma)\cap B_t(0)}m(\lambda)\le 
C t^{n-2}
\end{equation}
\cite[Theorem 8.3]{DKV}. The estimations (\ref{2.18}) and (\ref{2.19}) 
contain more information about the distribution of $\Lambda(\Gamma)$ then 
just the Weyl law. Indeed,  the
eigenvalue of $\Delta$ corresponding to $\lambda\in\Lambda_{\temp}(\Gamma)$ 
equals $\parallel\lambda\parallel^2+\parallel\rho\parallel^2$. So if we choose
$\Omega$ in (\ref{2.18}) to be the unit ball, then (\ref{2.18}) together with
(\ref{2.19}) reduces to Weyl's law for $\Gamma\bs S$. 

We note that  (\ref{2.18}) and (\ref{2.19}) can also be rephrased in terms of
representation theory. Let $R$ be the right regular representation of
$G$ in $L^2(\Gamma\bs G)$ defined by
\[(R(g_1)f)(g_2)=f(g_2g_1),\quad f\in L^2(\Gamma\bs G),\;g_1,g_2\in G.\]
Let $\hat G$ be the unitary dual of $G$, i.e., the set of equivalence classes
of irreducible unitary representations of $G$.
Since $\Gamma\bs G$ is compact, it is well known that $R$ decomposes into 
direct sum of irreducible unitary representations of $G$. Given $\pi\in\hat G$,
let $m(\pi)$ be the multiplicity with which $\pi$ occurs in $R$. Let $\H_\pi$
denote the Hilbert space in which $\pi$ acts. Then
\[L^2(\Gamma\bs G)\cong\bigoplus_{\pi\in\hat G} m(\lambda)\H_\pi.\]
Now observe that $L^2(\Gamma\bs S)=L^2(\Gamma\bs G)^K$. Let $\H_\pi^K$ denote 
the subspace of $K$-fixed vectors in $\H_\pi$. Then
\[L^2(\Gamma\bs S)\cong \bigoplus_{\pi\in\hat G} m(\lambda)\H_\pi^K.\]
Note that $\dim\H^K_\pi\le 1$. Let $\hat G(1)\subset \hat G$ be the subset of 
all $\pi$ with $\H^K_\pi\not=\{0\}$. This is the spherical dual. Given $\pi\in
\hat G$, let $\lambda_\pi$ be the infinitesimal character of $\pi$. If $\pi\in
\hat G(1)$, then $\lambda_\pi\in \af^\ast_\C/W$. Moreover $\pi\in\hat G(1)$ is
tempered, if $\pi$ is unitarily induced from the minimal parabolic subgroup
$P=NA$. In this case we have $\lambda_\pi\in i\af^\ast/W$. So (\ref{2.18}) can 
be rewritten as
\begin{equation}\label{2.20}
\sum_{\substack{\pi\in\hat G(1)\\
\lambda_\pi\in it\Omega}}m(\pi)=
\frac{\vol(\Gamma\bs G)}{|W|}\int_{it\Omega}\beta(\lambda)\,d\lambda
+O(t^{n-1}),\quad t\to\infty.
\end{equation}
\section{Automorphic forms}
\setcounter{equation}{0} 

The theory of automorphic forms is concerned with harmonic analysis on locally
symmetric spaces $\Gamma\bs S$ of finite volume. Of particular interest are
arithmetic groups $\Gamma$. This means that we consider a connected
semi-simple algebraic group $\G$ defined over $\Q$ such that $G=\G(\R)$ 
 and $\Gamma$ is a subgroup of $\G(\Q)$ which
is commensurable with $\G(\Z)$, where $\G(\Z)$ is defined with respect to some
embedding $\G\subset \GL(m)$. The standard example is $\G=\SL(n)$ and
$\Gamma(N)\subset \SL(n,\Z)$ the principal congruence subgroup of level $N$.
A basic feature of 
arithmetic groups is that the quotient $\Gamma\bs S$ has finite volume 
\cite{BH}. Moreover in many important cases it is non-compact. A typical 
example for that is  $\Gamma(N)\bs \SL(n,\R)/\SO(n)$.

In this section we discuss only the case of the upper half-plane $\bH$ and
we consider  congruence subgroups of $\SL(2,\Z)$. For
$N\ge 1$ the principal congruence subgroup $\Gamma(N)$ of level $N$ is defined
as
\[\Gamma(N)=\bigl\{\gamma\in\SL(2,\Z)\colon \gamma\equiv\Id\mod N\bigr\}.\]
A congruence subgroup $\Gamma$ of $\SL(2,\Z)$ is a subgroup for which there 
exists $N\in\N$ such that $\Gamma$ contains $\Gamma(N)$. An example of a 
congruence subgroup is the Hecke group
\[\Gamma_0(N)=\left\{\begin{pmatrix}a&b\\c&d\end{pmatrix}\in\SL(2,\Z)\colon
c\equiv 0\mod N\right\}.\] 
If $\Gamma$ is torsion free, the quotient $\Gamma\bs \bH$ is a finite area, 
non-compact, hyperbolic surface. It has a decomposition 
\begin{equation}\label{3.1}
\Gamma\bs \bH=M_0\cup Y_1\cup\cdots\cup Y_m,
\end{equation}
into the union of a compact surface with boundary $M_0$ and a finite number 
of ends 
$Y_i\cong [a,\infty)\times S^1$
which are equipped with the Poincar\'e metric. In general, $\Gamma\bs \bH$ may
have a finite number of quotient singularities. The quotient $\Gamma(N)\bs\bH$
is the modular surface  $X(N)$. 

Let $\Delta$ be the hyperbolic
Laplacian (\ref{2.1b}). A Maass automorphic form is a smooth function
$f\colon \bH\to\C$ which satisfies
\begin{enumerate}
\item[a)] $f(\gamma z)=f(z)$, $\gamma\in\Gamma$.
\item[b)] There exists $\lambda\in\C$ such that $\Delta f=\lambda f$.
\item[c)] $f$ is slowly increasing. 
\end{enumerate}
Here the last condition means that there exist $C>0$ and $N\in\N$ such that the
restriction $f_i$ of $f$ to $Y_i$ satisfies
\[|f_i(y,x)|\le Cy^N,\quad y\ge a,\; i=1,...,m.\]
Examples are the Eisenstein series. Let $a_1,...,a_m\in\R\cup\{\infty\}$ be
representatives of the $\Gamma$-conjugacy classes of parabolic fixed points of
$\Gamma$. The $a_i$'s are called {\it cusps}. For each $a_i$ let $\Gamma_{a_i}$
be the stabilizer of $a_i$ in $\Gamma$. Choose 
$\sigma_i\in\SL(2,\R)$ such that  
\[\sigma_i(\infty)=a_i,\quad \sigma_i^{-1}\Gamma_{a_i}\sigma_i=\left\{
\begin{pmatrix}1&n\\0&1\end{pmatrix}\colon n\in\Z\right\}.\]
Then the Eisenstein series $E_i(z,s)$ associated to the cusp $a_i$ is defined as
\begin{equation}
E_i(z,s)=\sum_{\gamma\in\Gamma_{a_i}\bs \Gamma}\Im(\sigma_i^{-1}\gamma z)^s,\quad 
\Re(s)>1.
\end{equation}
The series converges absolutely and uniformly on compact subsets of the 
half-plane $\Re(s)>1$ and it satisfies the following properties.
\begin{enumerate}
\item[1)] $E_i(\gamma z,s)=E_i(z,s)$ for all $\gamma\in\Gamma$.
\item[2)] As a function of $s$, $E_i(z,s)$ admits a meromorphic continuation 
to $\C$ which is regular on the line $\Re(s)=1/2$.
\item[3)] $E_i(z,s)$ is a smooth function of $z$ and satisfies
$\Delta_zE_i(z,s)=s(1-s) E_i(z,s).$
\end{enumerate}
As example consider the modular group $\Gamma(1)$ which has a single cusp 
$\infty$. The Eisenstein series attached to $\infty$ is the well-known 
series
\[E(z,s)=\sum_{\substack{(m,n)\in\Z^2\\(m,n)=1}}\frac{y^s}{|mz+n|^{2s}}.\]
In the general case, the Eisenstein series were first studied by Selberg 
\cite{Se1}. The Eisenstein series are closely related with the study of
the spectral resolution of $\Delta$.  Regarded
as unbounded operator 
\[\Delta\colon C_c^\infty(\Gamma\bs \bH)\to L^2(\Gamma\bs \bH),\]
$\Delta$ is essentially self-adjoint \cite{Roe}. Let $\bar\Delta$ be
the  unique self-adjoint extension  of $\Delta$. The important new feature due
to the non-compactness of $\Gamma\bs\bH$ is that $\bar\Delta$ has a large
continuous spectrum which is governed by the Eisenstein series. The following
basic result is due to Roelcke \cite{Roe}.
\begin{prop}\label{p3.1}
The spectrum of $\bar\Delta$ is the union of a pure point spectrum 
$\sigma_{pp}(\bar\Delta)$ and an absolutely continuous spectrum
$\sigma_{ac}(\bar\Delta)$.\\
1) The pure point spectrum consists of eigenvalues 
$0=\lambda_0<\lambda_1\le\cdots$ of finite multiplicities with no finite points 
of accumulation.\\
2) The absolutely continuous spectrum equals $[1/4,\infty)$ with
multiplicity equal to the number of cusps of $\Gamma\bs \bH$.
\end{prop}
Of particular interest are the eigenfunctions of $\bar\Delta$. They are 
Maass automorphic forms. This can be seen by studying the Fourier expansion 
of an eigenfunction in the cusps. As an example consider 
$f\in C^\infty(\Gamma_0(N)\bs \bH)$ which satisfies
\[\Delta f=\lambda f, \quad f(z)=f(-\bar z),\quad \int_{\Gamma_0(N)\bs \bH}
|f(z)|^2\;dA(z)<\infty.\]
Assume that $\lambda=1/4+r^2$, $r\in\R$. Then $f(x+iy)$ admits  a Fourier 
expansion 
w.r.t. $x$ of the form
\begin{equation}\label{3.2}
f(x+iy)=\sum_{n=1}^\infty a(n)\sqrt{y}K_{ir}(2\pi ny)\cos(2\pi nx),
\end{equation}
where $K_\nu(y)$ is the modified Bessel function which may be defined by
\[K_\nu(y)=\int_0^\infty e^{-y\cosh t}\cosh(\nu t)\;dt\]
and it satisfies
\[K_\nu^{\prime\prime}(y)+\frac{1}{y}K_\nu^\prime(y)+
\left(1-\frac{\nu^2}{y^2}\right)K_\nu(y)=0.\]
Now note that $K_\nu(y)=O(e^{-cy})$ as $y\to\infty$. This implies that $f$ 
is rapidly decreasing in the  cusp $\infty$. A similar Fourier expansion 
holds in the other cusps. This implies that $f$ is rapidly decreasing in all 
cusps and therefore, it is a Maass automorphic form. In fact, since the zero 
Fourier coefficients vanish in all cusps, $f$ is a Maass {\it cusp form}.
In general, the space of cusp forms $L^2_{\cu}(\Gamma\bs\bH)$ is defined as
the subspace of all $f\in L^2(\Gamma\bs\bH)$ such that for almost all 
$y\in\R^+$:
\[\int_0^1 f(\sigma_k(x+iy))\;dx=0,\quad k=1,...,m.\]
This is an invariant subspace of $\bar\Delta$ and the restriction of 
$\bar\Delta$ to $L^2_{\cu}(\Gamma\bs \bH)$ has pure point spectrum, i.e., 
$L^2_{\cu}(\Gamma\bs \bH)$ is the span of square integrable eigenfunctions of 
$\Delta$. Let $L^2_{\res}(\Gamma\bs\bH)$ be the orthogonal complement of
$L^2_{\cu}(\Gamma\bs\bH)$ in $L^2(\Gamma\bs\bH)$. Thus
\[L^2_{pp}(\Gamma\bs\bH)=L^2_{\cu}(\Gamma\bs\bH)\oplus 
L^2_{\res}(\Gamma\bs\bH).\]
The subspace $L^2_{\res}(\Gamma\bs\bH)$ can be described as follows. The poles
of the Eisenstein series $E_i(z,s)$ in the half-plane $\Re(s)>1/2$ are all
simple and are contained in the interval $(1/2,1]$. Let $s_0\in(1/2,1]$ be a
pole of $E_i(z,s)$ and put
\[\psi=\mathrm{Res}_{s=s_0}E_i(z,s).\]
Then $\psi$ is a square integrable eigenfunction of $\Delta$ with eigenvalue
$\lambda=s_0(1-s_0)$. The set of all such residues of the Eisenstein series
$E_i(z,s)$, $i=1,...,m$, spans $L^2_{\res}(\Gamma\bs\bH)$. This is a 
finite-dimensional space which  is called the {\it residual subspace}. 
The corresponding eigenvalues form the {\it residual spectrum} of 
$\bar\Delta$. 
So we are left with the cuspidal eigenfunctions or {\it Maass cusp forms}.
Cusp forms are the building blocks of the theory of automorphic forms. They
play an important role in number theory. To illustrate this consider an even 
Maass cusp form $f$ for $\Gamma(1)$ with eigenvalue $\lambda=1/4+r^2$, $r\in\R$.
Let $a(n)$, $n\in\N$, be the Fourier coefficients of $f$ given by (\ref{3.2}). 
Put
\[L(s,f)=\sum_{n=1}^\infty\frac{a(n)}{n^s},\quad\Re(s)>1.\]
This Dirichlet series converges absolutely and uniformly in the half-plane
$\Re(s)>1$. Let
\begin{equation}\label{3.2a}
\Lambda(s,f)=\pi^{-s}{\bf \Gamma}\left(\frac{s+ir}{2}\right)
{\bf \Gamma}\left(\frac{s-ir}{2}\right) L(s,f).
\end{equation}
Then the modularity of $f$ implies that $\Lambda(s,f)$ has an analytic 
continuation to the whole complex plane and satisfies the functional equation
\[\Lambda(s,f)=\Lambda(1-s,f)\]
\cite[Proposition 1.9.1]{Bu}. 
Under additional assumptions on $f$, the Dirichlet series
$L(s,f)$ is also an Euler product. This is related to the arithmetic nature 
of the groups $\Gamma(N)$. The surfaces $X(N)$ carry a family of
algebraically defined operators $T_n$, the so called {\it Hecke operators},
which  for $(n,N)=1$ are defined by
\[T_nf(z)=\frac{1}{\sqrt{n}}\sum_{\substack{ad=n\\b\mod d}} 
f\left(\frac{az+b}{d}\right).\]
These are closely related to the cosets of the finite index subgroups
\[\begin{pmatrix}n&0\\0&1\end{pmatrix}\Gamma(N)
\begin{pmatrix}n&0\\0&1\end{pmatrix}^{-1}\cap \Gamma(N)\]
of $\Gamma(N)$. Each  $T_n$ defines a linear transformation of $L^2(X(N))$.
The $T_n$, $n\in\N$,  are a commuting family of normal operators which 
also commute with $\Delta$. Therefore, each $T_n$ leaves the eigenspaces of
$\Delta$ invariant. So we may assume that $f$ is a common eigenfunction of
$\Delta$ and $T_n$, $n\in\N$:
\[\Delta f=(1/4+r^2)f,\quad T_n f=\lambda(n) f.\] 
If $f\not=0$, then $a(1)\not=0$. So we can normalize $f$ such that $a(1)=1$. 
Then it follows that $a(n)=\lambda(n)$ and 
the Fourier coefficients satisfy the following  multiplicative relations
\[a(m)a(n)=\sum_{d|(m,n)}a\left(\frac{mn}{d^2}\right).\]
This implies that  $L(s,f)$ is an Euler product
\begin{equation}
L(s,f)=\sum_{n=1}^\infty a(n) n^{-s}
=\prod_p\left(1-a(p) p^{-s}+p^{-2s}\right)^{-1},
\end{equation}
which converges for $\Re(s)>1$. $L(s,f)$ is the basic example of an automorphic
$L$-function.  It is convenient to write this Euler product in a different way.
Introduce roots $\alpha_p,\beta_p$ by
\[\alpha_p\beta_p=1,\quad \alpha_p+\beta_p=a(p).\]
Let 
\[A_p=\begin{pmatrix}\alpha_p&0\\0&\beta_p\end{pmatrix}.\]
Then
\[L(s,f)=\prod_p\det\left(\Id-A_pp^{-s}\right)^{-1}.\]
Now let $\rho\colon \GL(2,\C)\to\GL(N,\C)$
be a representation. Then we can form a new Euler product by
\[L(s,f,\rho)=\prod_p\det\left(\Id-\rho(A_p)p^{-s}\right)^{-1},\]
which converges in some half-plane.
It is part of the general conjectures of Langlands \cite{La2} that each of 
these Euler products admits a meromorphic extension to $\C$ and satisfies a
functional equation. The construction of Euler products for Maass
cusp forms can be extended to other groups $\G$, in particular to cusp forms
on $\GL(n)$. 
It is also conjectured that $L(s,f,\rho)$ is an
automorphic $L$-function of an automorphic form on some $\GL(n)$. This is
part of the functoriality principle of Langlands. Furthermore, all 
$L$-functions that occur in number theory and algebraic geometry are expected
to be automorphic $L$-functions.  
This is one of the main reasons for the interest in the study of cusp forms.
Other applications are discussed in \cite{Sa1}. 

\section{The Weyl law and existence of cusp forms}
\setcounter{equation}{0} 

Since $\Gamma(N)\bs\bH$ is not compact, it is not clear that there exist any
eigenvalues $\lambda>0$.  
By Proposition \ref{p3.1} the continuous spectrum of $\bar\Delta$ equals
$[1/4,\infty)$. Thus all eigenvalues $\lambda\ge 1/4$ are embedded in the
continuous spectrum. It is well-known in mathematical physics, that
embedded eigenvalues are unstable under perturbations and therefore, are
difficult to study. 

One of the basic tools to study the cuspidal spectrum is the Selberg trace
formula \cite{Se2}. The new terms in the trace formula, which are due to
the non-compactness of $\Gamma\bs\bH$ arise from the parabolic conjugacy
classes in $\Gamma$ and the Eisenstein series. The contribution of the
Eisenstein series is given by their zeroth Fourier coefficients of the
Fourier expansion in the cusps. The 
zeroth Fourier coefficient of the Eisenstein series $E_k(z,s)$ in the cusp 
$a_l$ is given by
\[\int_0^1 E_k(\sigma_l(x+iy),s)\;dx=y^s+C_{kl}(s)y^{1-s},\]
where $C_{kl}(s)$ is a meromorphic function of $s\in\C$. Put
\[C(s):=\left(C_{kl}(s)\right)_{k,l=1}^m.\]
This is the so called scattering matrix. Let 
\[\phi(s):=\det C(s).\]
Let the notation be as in (\ref{2.2}) and assume that $\Gamma$ has no torsion. 
Then the trace formula is the following identity.
\begin{equation}\label{4.1}
\begin{split}
\sum_{j} h(r_j)&=\frac{\Area(\Gamma\bs \bH)}{4\pi}\int_\R h(r)r
\tanh(\pi r)\;dr
+\sum_{\{\gamma\}_\Gamma} 
\frac{l(\gamma_0)}{2\sinh \left( \frac{l(\gamma)}{2}\right)} g(l(\gamma))\\
&\quad +\frac{1}{4\pi} \int^\infty_{-\infty}h(r)\frac{\phi^\prime}
{\phi}(1/2+ir)\;dr  
-\frac{1}{4}\phi(1/2)h(0)\\
& \quad -\frac{m}{2\pi}\int^\infty_{-\infty} h(r)
\frac{{\bf \Gamma}^\prime}{{\bf \Gamma}}(1+ir)dr +\frac{m}{4} h(0)-m\ln 2\; g(0).
\end{split}
\end{equation}
The trace formula holds for every discrete subgroup $\Gamma\subset \SL(2,\R)$ 
with finite coarea. 
In analogy to the counting function of the eigenvalues we
introduce the winding number 
\[M_\Gamma(\lambda)=-\frac{1}{4\pi}\int_{-\lambda}^\lambda\frac{\phi^\prime}
{\phi}(1/2+ir)\;dr\]
which measures the continuous spectrum. 
Using the cut-off Laplacian of 
Lax-Phillips \cite{CV} one can deduce the following elementary bounds 
\begin{equation}\label{4.1a}
N_\Gamma(\lambda)\ll\lambda^2,\quad M_\Gamma(\lambda)\ll\lambda^2,\quad 
\lambda\ge 1.
\end{equation}
These bounds imply that the trace formula (\ref{4.1}) holds for a larger
class of functions. In particular, it can be applied to the heat 
kernel $k_t(u)$.  Its spherical Fourier transform 
equals $h_t(r)=e^{-t(1/4+r^2)}$, $t>0$.
If we insert $h_t$ into the trace formula we get the following asymptotic 
expansion as $t\to 0$.
\begin{equation}\label{4.3}
\begin{split}
\sum_{j}e^{-t\lambda_j}-\frac{1}{4\pi}\int_\R &e^{-t(1/4+r^2)}
\frac{\phi^\prime}{\phi}(1/2+ir)\;dr\\
&=\frac{\Area(\Gamma\bs\bH)}{4\pi t}+\frac{a\log t}{\sqrt{t}}+
\frac{b}{\sqrt{t}}+O(1)
\end{split}
\end{equation} 
for certain constants $a,b\in\R$. Using \cite[(8.8), (8.9)]{Se2} it 
follows that the winding number $M_\Gamma(\lambda)$ is monotonic increasing 
for $r\gg0$. Therefore we can apply a Tauberian theorem to (\ref{4.3}) and we
get the  Weyl law (\ref{1.4}).

In general, we cannot estimate separately the counting function and the winding
number. For congruence subgroups, however, the entries of
the scattering matrix can be expressed in terms of well-known analytic
functions. For $\Gamma(N)$ the determinant of the scattering 
matrix $\phi(s)$ has been computed by Huxley \cite{Hu}. It has the 
form
\begin{equation}\label{4.5}
\phi(s)=(-1)^lA^{1-2s}\left(\frac{{\bf \Gamma}(1-s)}{{\bf \Gamma}(s)}\right)^k
\prod_\chi\frac{L(2-2s,\bar\chi)}{L(2s,\chi)},
\end{equation}
where $k,l\in\Z$, $A>0$, the product runs
over Dirichlet characters $\chi$ to some modulus dividing $N$ and $L(s,\chi)$
is the Dirichlet $L$-function with character $\chi$. Especially for $\Gamma(1)$
we have
\begin{equation}\label{4.5a}
\phi(s)=\sqrt{\pi}\frac{{\bf \Gamma}(s-1/2)\zeta(2s-1)}{{\bf \Gamma}(s)
\zeta(2s)},
\end{equation}
where $\zeta(s)$ denotes the Riemann zeta function.

Using Stirling's 
approximation formula to estimate the logarithmic derivative of the Gamma 
function and standard estimations for the logarithmic derivative
of Dirichlet $L$-functions on the line $\Re(s)=1$ \cite[Theormem 7.1]{Pr}, 
we get
\begin{equation}\label{4.6}
\frac{\phi'}{\phi}(1/2+ir)=O(\log(4+|r|)),\quad |r|\to\infty.
\end{equation}
This implies that
\begin{equation}\label{4.7}
M_{\Gamma(N)}(\lambda)\ll\lambda\log\lambda.
\end{equation}
Together with (\ref{1.4}) we obtain Weyl's law for the point spectrum
\begin{equation}\label{4.8}
N_{\Gamma(N)}(\lambda) \sim
\frac{\Area(X(N))}{4\pi}\lambda^2,\quad\lambda\to\infty,
\end{equation}
which is due to Selberg \cite[p.668]{Se2}.
A similar formula holds for other congruence groups such as $\Gamma_0(N)$.
In particular,  (\ref{4.8})  implies
that for congruence groups $\Gamma$ there exist infinitely many linearly 
independent Maass cusp forms.

A  proof of the Weyl law (\ref{4.8}) which avoids the use of the
constant terms of the Eisenstein series has recently been given by 
Lindenstrauss and Venkatesh \cite{LV}. Their approach is based on the 
construction of convolution operators with purely cuspidal image.

Neither of these methods  give good  estimates of the remainder term. 
One approach to obtain estimates of the remainder term is based on the Selberg 
zeta function
\[
Z_\Gamma(s)=\prod_{\{\gamma\}_\Gamma}\prod_{k=0}^\infty\left(1-e^{-(s+k)\ell(\gamma)}
\right),\quad \Re(s)>1,
\]
where the outer product runs over the primitive hyperbolic conjugacy classes
in $\Gamma$ and $\ell(\gamma)$ is the length of the closed
geodesic associated to $\{\gamma\}_\Gamma$. The infinite product converges  
absolutely in the indicated half-plane and admits an analytic continuation to
the whole complex plane. 
If $\lambda=1/4+r^2$, $r\in
\R\cup i(1/2,1]$, is an  eigenvalue of $\Delta$, then $s_0=1/2+ir$ is a zero of 
$Z_\Gamma(s)$. Using this fact and standard methods of analytic number theory
one can derive  the following strong form of the Weyl law
\cite[Theorem 2.28]{Hj}, \cite[Theorem 7.3]{Ve}.
\begin{theo}\label{th4.1} 
Let $m$ be the number of cusps of $\Gamma\bs\bH$. There exists
$c>0$ such that
\[N_\Gamma(\lambda)+M_\Gamma(\lambda)=\frac{\Area(\Gamma\bs \bH)}{4\pi}\lambda^2
-\frac{m}{\pi}\lambda\log\lambda+c\lambda+
O\left(\lambda(\log\lambda)^{-1}\right)\]
as $\lambda\to\infty$.
\end{theo}

Together with (\ref{4.7})  we obtain Weyl's law with remainder term.
\begin{theo}\label{th4.2} 
For every $N\in\N$ we have
\[N_{\Gamma(N)}(\lambda)=\frac{\Area(X(N))}{4\pi}\lambda^2+
O(\lambda\log\lambda)\]
as $\lambda\to\infty$.
\end{theo}

The  use of the Selberg zeta function to estimate the remainder term
is  limited to  rank one cases. However,
the remainder term  can also be estimated by H\"ormander's
method using the trace formula as in the compact case. We describe the main
steps.
Let $\varepsilon>0$ such
that $\ell(\gamma)>\varepsilon$ for all hyperbolic conjugacy classes
$\{\gamma\}_{\Gamma(N)}$.
 Choose $g\in C_c^\infty(\R)$ to be even and such that $\supp g\subset 
(-\varepsilon,\varepsilon)$. Let $h(z)=\int_\R g(r)e^{-irz}\,dr$. Then the
hyperbolic contribution in the trace formula (\ref{4.1}) drops out. We
symmetrize the spectrum by $r_{-j}=-r_j$, $j\in\N$. Then for each $t\in\R$ we 
have
\begin{equation}\label{4.8a}
\begin{split}
\sum_{j} h(t-r_j)&=\frac{\Area(X(N))}{2\pi}\int_\R h(t-r)r
\tanh(\pi r)\;dr\\
& +\frac{1}{2\pi} \int^\infty_{-\infty}h(t-r)\frac{\phi^\prime}
{\phi}(1/2+ir)\;dr  
-\frac{1}{2}\phi(1/2)h(t)\\
&  -\frac{m}{\pi}\int^\infty_{-\infty} h(t-r)
\frac{{\bf \Gamma}^\prime}{{\bf \Gamma}}(1+ir)dr +\frac{m}{2} 
h(t)-2m\ln 2\; g(0).
\end{split}
\end{equation}
Now we need to estimate  the  behavior of the 
terms on the right hand side as $|t|\to\infty$. The first integral has been
already considered in (\ref{2.5}). It is of order $O(|t|)$.  
To deal with the second integral we use (\ref{4.6}).
This implies
\begin{equation}\label{4.8b}
\int_\R h(t-r)\frac{\phi'}{\phi}(1/2+ir)\;dr=O(\log(|t|)),\quad
|t|\to\infty. 
\end{equation}
Using Stirling's formula we get
\begin{equation*}
\int_\R h(t-r)\frac{\Gamma'}{\Gamma}(1+ir)\;dr=O(\log(|t|)),\quad
|t|\to\infty. 
\end{equation*}
The remaining terms are bounded. 
Combining these estimations, we get
\[\sum_{j}h(t-r_j)=O(|t|),\quad |t|\to\infty.\] 
Therefore,  Lemma \ref{l2.2} holds also in the present case. It remains to
establish the analog of Lemma \ref{l2.4}. Using (\ref{2.9}) and (\ref{4.6})
we get
\begin{equation}\label{4.9}
\int_{-\lambda}^\lambda\int_\R h(t-r)\frac{\phi'}{\phi}(1/2+it)\;dr\;dt=
O(\lambda\log\lambda).
\end{equation}
Similarly, using Stirling's  formula and (\ref{2.9}), we obtain
\begin{equation}\label{4.10}
\int_{-\lambda}^\lambda\int_\R 
h(t-r)\frac{{\bf \Gamma}'}{{\bf \Gamma}}(1+it)\;dr\;dt=
O(\lambda\log\lambda).
\end{equation}
The integral of the remaining terms is of order $O(\lambda)$. Thus we obtain
\begin{equation}\label{4.11}
\int_{-\lambda}^\lambda \sum_{j=-\infty}^\infty h(t-r_j)\;dt=
\frac{\Area(X(N))}{2\pi}\lambda^2+O(\lambda\log \lambda)
\end{equation}
as $\lambda\to\infty$. Now we proceed in exactly the same way as in the
compact case. Using Lemma \ref{l2.3} and (\ref{4.11}), Theorem \ref{th4.2}
follows.
\hfill$\Box$

The Weyl law shows that for congruence groups Maass cusp forms exist in 
abundance. In general very little is
known. Let $\Gamma$ be any discrete, co-finite subgroup of $\SL(2,\R)$. Then
by Donnelly \cite{Do}  the following general bound is known 
\[\limsup_{\lambda\to\infty}\frac{N_\Gamma^{\cu}(\lambda)}{\lambda}\le 
\frac{\Area(\Gamma\bs\bH)}
{4\pi}.\]
A group $\Gamma$ for which the equality is attained
is called {\it essentially cuspidal} by Sarnak \cite{Sa2}. 
By (\ref{4.8}),  $\Gamma(N)$ is essentially cuspidal. 
The study of the behavior of eigenvalues under deformations of $\Gamma$,
initiated Phillips and Sarnak \cite{PS1}, \cite{PS2}, supports the 
conjecture that essential cuspidality may be limited to special arithmetic
groups. 

The consideration of the behavior of cuspidal eigenvalues under
deformations was started by Colin de Verdiere \cite{CV} in the more general 
context of metric perturbations. One of his main results \cite[Th\'eor\`eme 7]
{CV} states that under a generic compactly supported conformal
perturbation of the hyperbolic metric of $\Gamma\bs\bH$ all Maass cusp forms
are dissolved. This means that each point $s_j=1/2+ir_j$, $r_j\in\R$, such that 
$\lambda_j=s_j(1-s_j)$ is an eigenvalue moves under the perturbation into 
the half-plane $\Re(s)<1/2$ and becomes a pole of the scattering matrix 
$C(s)$. 

In the present context we are only interested in deformations such that the 
curvature stays constant. Such deformations are given by curves in 
 the Teichm\"uller space $\cT(\Gamma)$  of $\Gamma$. The space
$\cT(\Gamma)$ is known to be a finite-dimensional  and therefore, it is 
by no means clear that the results of \cite{CV} will continue to hold for
perturbations of this restricted type. For $\Gamma(N)$ this problem
has been studied in \cite{PS1}, \cite{PS2}. One of the main results is an 
analog of Fermi's golden rule which gives a sufficient condition for 
a cusp form of $\Gamma(N)$ to be dissolved under a deformation in 
$\cT(\Gamma(N))$. Based on these results, Sarnak made the following conjecture
\cite{Sa2}:

\noindent
{\bf Conjecture.}
\begin{enumerate}
\item[(a)] {\it The generic $\Gamma$ in a given Teichm\"uller space of finite 
area hyperbolic surfaces is not essentially cuspidal.}
\item[(b)] {\it Except for the Teichm\"uller space of the once punctured torus,
the generic $\Gamma$ has only finitely many  eigenvalues. }
\end{enumerate}

\section{Higher rank}
\setcounter{equation}{0} 

In this section we consider an arbitrary locally symmetric space $\Gamma\bs S$
defined by an arithmetic subgroup $\Gamma\subset \G(\Q)$, where  
$\G$ is a semi-simple algebraic group over $\Q$ with finite center, 
$G=\G(\R)$ and $S=G/K$. 
The basic example will be $\G=\SL(n)$ and $\Gamma=\Gamma(N)$, the principal
congruence subgroup of $\SL(n,\Z)$ of level $N$ which consists of all $\gamma
\in \SL(n,\Z)$ such that $\gamma\equiv \Id\mod N$. 

Let $\Delta$ be the Laplacian of $\Gamma\bs S$, and let $\bar\Delta$ be the 
closure of $\Delta$ in $L^2$. Then $\bar\Delta$ is a non-negative self-adjoint
operator in $L^2(\Gamma\bs S)$. The properties of its spectral resolution can 
be derived from the known structure of the spectral resolution of the regular
representation $R_\Gamma$ of $G$ on $L^2(\Gamma\bs G)$ \cite{La1}, \cite{BG}.
In this way we get the following generalization of Proposition \ref{p3.1}.
\begin{prop}\label{p5.1}
The spectrum of $\bar\Delta$ is the union of a point spectrum 
$\sigma_{pp}(\bar\Delta)$ and an absolutely continuous spectrum
$\sigma_{ac}(\bar\Delta)$.\\
1) The point spectrum consists of eigenvalues 
$0=\lambda_0<\lambda_1\le\cdots$ of finite multiplicities with no finite
point of accumulation.\\
2) The absolutely continuous spectrum equals $[b,\infty)$ for some $b>0$.
\end{prop}
The theory of Eisenstein series \cite{La1}
provides a complete set of generalized eigenfunctions for $\Delta$. The
corresponding wave packets span the absolutely continuous subspace 
$L^2_{\ac}(\Gamma\bs S)$.  This allows us to determine the constant $b$ 
explicitly in terms of the root structure. The statement about the point
spectrum was proved in \cite[Theorem 5.5]{BG}.

Let $L^2_{\di}(\Gamma\bs S)$ be the closure of the span of all eigenfunctions.
It contains the subspace of cusp forms $L^2_{\cu}(\Gamma\bs S)$. We recall
its definition. Let $\bP\subset \G$ be a
parabolic subgroup defined over $\Q$ \cite{Bo}. Let $P=\bP(\R)$. This is
a cuspidal parabolic subgroup of $G$ and all cuspidal parabolic subgroups 
arise in this way. Let ${\bf N}_P$ be the unipotent radical of $\bP$ and let
$N_P={\bf N}_P(\R)$.
Then $N_P\cap\Gamma\bs N_P$ is compact. A cusp form is a smooth function $\phi$
on $\Gamma\bs S$ which is a joint eigenfunction of the ring $\cD(S)$ of 
invariant differential operators on $S$, and which satisfies 
\begin{equation}\label{5.1}
\int_{N_P\cap\Gamma\bs N_P} \phi(nx)\,dn=0
\end{equation}
for all  cuspidal parabolic subgroups $P\not=G$. Each cusp form is rapidly
decreasing and hence 
square integrable. Let $L^2_{\cu}(\Gamma\bs S)$ be the closure in 
$L^2(\Gamma\bs S)$ of the linear span of all cusp forms. Then 
$L^2_{\cu}(\Gamma\bs S)$ is an invariant subspace of $\bar\Delta$ which is
contained in $L^2_{\di}(\Gamma\bs S)$. 

Let $L^2_{\res}(\Gamma\bs S)$ be the orthogonal complement of 
$L^2_{\cu}(\Gamma\bs S)$ in  $L^2_{\di}(\Gamma\bs S)$, i.e., we have an 
orthogonal decomposition
\[L^2_{\di}(\Gamma\bs S)=L^2_{\cu}(\Gamma\bs S)\oplus L^2_{\res}(\Gamma\bs S).\]
It follows from Langlands's theory of Eisenstein systems that
$L^2_{\res}(\Gamma\bs S)$ is spanned by {\it iterated residues} of cuspidal 
Eisenstein series \cite[Chapter 7]{La1}. Therefore $L^2_{\res}(\Gamma\bs S)$ 
is called the {\it residual subspace}. 

Let $N_\Gamma^{\di}(\lambda)$, $N_\Gamma^{\cu}(\lambda)$, and 
$N_\Gamma^{\res}(\lambda)$ be the 
counting function of the eigenvalues with eigenfunctions belonging to the
corresponding subspaces. The following general results about the growth of
the counting functions are known for any lattice $\Gamma$ in a real
semi-simple Lie group. Let $n=\dim S$. Donnelly \cite{Do} has proved the
following bound for the cuspidal spectrum
\begin{equation}\label{5.2}
\limsup_{\lambda\to\infty}\frac{N_\Gamma^{\cu}(\lambda)}{\lambda^{n}}\le
\frac{\vol(\Gamma\bs S)}{(4\pi)^{n/2}{\bf \Gamma}\left(\frac{n}{2}+1\right)},
\end{equation}
where ${\bf \Gamma}(s)$ denotes the Gamma function. For  the full discrete 
spectrum, we have at least an upper bound for the
growth of the counting function. The main result of \cite{Mu2} states that 
\begin{equation}\label{5.3}
N_\Gamma^{\di}(\lambda)\ll(1+\lambda^{4n}).
\end{equation}
This result implies that invariant integral operators are trace class on
the discrete subspace which is the starting point for the trace formula.
 The proof of (\ref{5.3}) relies on the description
of the residual subspace in terms of iterated residues of Eisenstein series.
One actually expects that the growth of the residual spectrum is of lower order
than the cuspidal spectrum. For $\SL(n)$ the residual spectrum has been 
determined by Moeglin and Waldspurger \cite{MW}. Combined with (\ref{5.2}) it
follows that for $G=\SL(n)$ we have
\begin{equation}\label{5.4}
N_{\Gamma(N)}^{\res}(\lambda)\ll \lambda^{d-1},
\end{equation}
where $d=\dim \SL(n,\R)/\SO(n)$. 

In \cite{Sa2} Sarnak  conjectured that if $\rk(G/K)>1$,
each  irreducible lattice $\Gamma$ in $G$ is essentially cuspidal in the
sense that Weyl's law holds for $N_\Gamma^{\cu}(\lambda)$, i.e., equality
holds in (\ref{5.2}). This conjecture
has now been established in quite generality. A. Reznikov proved it for 
congruence
groups in a group $G$ of real rank one, S. Miller \cite{Mi} proved  it
for $\G=\SL(3)$ and $\Gamma=\SL(3,\Z)$, the author \cite{Mu3} established it 
for $\G=\SL(n)$ and a congruence group $\Gamma$. The method of \cite{Mu3}
is an extension of the heat equation method described in the previous 
section for the case of the upper half-plane. More recently, Lindenstrauss
and Venkatesh \cite{LV} proved the following result.
\begin{theo}
Let $\G$ be a split adjoint semi-simple group over $\Q$ and let $\Gamma\subset
\G(\Q)$ be a congruence subgroup. Let $n=\dim S$. Then 
\[N_\Gamma^{\cu}(\lambda)\sim \frac{\vol(\Gamma\bs S)}{(4\pi)^{n/2}
{\bf \Gamma}\left(\frac{n}{2}+1\right)}\lambda^{n},\quad \lambda\to\infty.\]
\end{theo}
The method is based on the construction of convolution operators with pure
cuspidal image. It avoids the delicate estimates of the contributions of the
Eisenstein series to the trace formula.
This proves existence of many cusp forms for these groups.

The next problem is to estimate the remainder term. For $\G=\SL(n)$, this
problem has been studied by E. Lapid and the author in \cite{LM}. 
Actually, we consider
not only the cuspidal spectrum of the Laplacian, but the cuspidal spectrum 
of the whole algebra of invariant differential operators.  

As $\cD(S)$ preserves the space of cusp
forms, we can proceed as in the compact case and decompose 
$L^2_{\cu}(\Gamma\bs S)$ into joint eigenspaces of $\cD(S)$. Given $\lambda\in
\af^\ast_\C/W$, let
\[
\E_{\cu}(\lambda)=\left\{\varphi\in L^2_{\cu}(\Gamma\bs S)\colon
D\varphi=\chi_\lambda(D)\varphi, \;D\in\cD(S)\right\}
\]
be the associated eigenspace.
Each eigenspace is finite-dimensional. Let $m(\lambda)=\E_{\cu}(\lambda)$. 
Define the cuspidal spectrum 
$\Lambda_{\cu}(\Gamma)$ to be
\[\Lambda_{\cu}(\Gamma)=\{\lambda\in\af^\ast_\C/W\colon m(\lambda)>0\}.\]
Then we have an orthogonal direct sum decomposition
\[L^2_{\cu}(\Gamma\bs S)=\bigoplus_{\lambda\in\Lambda_{\cu}(\Gamma)}
\E_{\cu}(\lambda).\]
Let the notation be as in ({2.18}) and ({2.19}).
Then in \cite{LM}  we established the following extension of 
main results of \cite{DKV} to congruence quotients of $S=\SL(n,\R)/\SO(n)$.
\begin{theo}\label{th5.3}
Let  $d=\dim S$. 
Let $\Omega\subset \af^\ast$ be a bounded domain with piecewise smooth
boundary. Then for $N\ge 3$ we have
\begin{equation}\label{5.5} 
\sum_{\lambda\in\Lambda_{\cu}(\Gamma(N)),
\lambda\in i t\Omega}m(\lambda)=
\frac{\vol(\Gamma(N)\bs S)}{|W|}\int_{t\Omega}\beta(i\lambda)\ d\lambda+
O\left(t^{d-1}(\log t)^{\max(n,3)}\right), 
\end{equation}
as $t\to\infty$, and
\begin{equation}\label{5.6} 
\sum_{\substack{\lambda\in\Lambda_{\cu}(\Gamma(N))\\
\lambda\in B_t(0)\setminus i\af^\ast}}m(\lambda)=
O\left(t^{d-2}\right),\quad t\to\infty.
\end{equation}
\end{theo}  
If we apply (\ref{5.5}) and (\ref{5.6}) to the unit ball in $\af^\ast$, we
get the following corollary.

\begin{corollary} Let $\G=\SL(n)$ and let $\Gamma(N)$ be the principal 
congruence
subgroup of $\SL(n,\Z)$ of level $N$. Let $S=\SL(n,\R)/\SO(n)$ and $d=\dim S$.
 Then for $N\ge 3$ we have
\[
N^{\cu}_{\Gamma(N)}(\lambda)=\frac{\vol(\Gamma(N)\bs S)}
{(4\pi)^{d/2}{\bf\Gamma}\left(\frac{d}{2}+1\right)}\lambda^d+
O\left(\lambda^{d-1}(\log \lambda)^{\max(n,3)}\right),\quad \lambda\to\infty.
\]
\end{corollary} 
The condition $N\ge 3$ is imposed for technical reasons. It guarantees that
the principal congruence subgroup $\Gamma(N)$ is neat in the sense of Borel,
and in particular, has no torsion. This simplifies the analysis by eliminating
the contributions of the non-unipotent conjugacy classes in the trace formula.

Note that $\Lambda_{\cu}(\Gamma(N))\cap i\af^\ast$ is the cuspidal tempered
spherical spectrum. The Ramanujan conjecture \cite{Sa3}
 for $\GL(n)$ at the Archimedean place states that
\[\Lambda_{\cu}(\Gamma(N))\subset i\af^\ast\]
so that (\ref{5.6}) is empty, if the Ramanujan conjecture is true. However, 
the Ramanujan conjecture is far from 
being proved. Moreover, it is known to be false for other groups $\G$ and
(\ref{5.6}) is what one can expect in general.

The method to prove Theorem \ref{th5.3} is an extension of the method of
\cite{DKV}. The Selberg trace formula, which is one of the basic tools in
 \cite{DKV}, needs to be replaced by the Arthur trace formula \cite{A1}, 
\cite{A2}. This requires to change the framework and to work with
 the adelic setting. It is also convenient to replace $\SL(n)$ by $\GL(n)$.

Again, one of the main issues is to estimate the terms in the trace formula
which are associated to Eisenstein series. Roughly speaking, these terms are 
a sum of integrals which generalize the integral
\[\int^\infty_{-\infty}h(r)\frac{\phi^\prime}{\phi}(1/2+ir)\;dr\]
in (\ref{4.1}). The sum is running over Levi components of parabolic subgroups
and square integrable automorphic forms on a given Levi component. 
The  functions which generalize $\phi(s)$ are obtained from the constant terms
of Eisenstein series. In general, they are difficult to describe. The main
ingredients are  logarithmic derivatives of
automorphic $L$-functions associated to automorphic forms on the Levi 
components. As example consider $\G=\SL(3)$, $\Gamma=\SL(3,\Z)$, and a 
standard maximal parabolic subgroup $P$ which has the form
\begin{equation*}
P=\left\{ \begin{pmatrix} m_1 & X\\
0 & m_2\end{pmatrix} \; \Big| \; m_i\in\GL(n_{i},\R),\;
\det m_1\cdot \det m_2=1\right\},
\end{equation*}
with $n_1+ n_2=3$. Thus there are exactly two standard maximal parabolic
subgroups. The standard Levi component of $P$ is 
\begin{equation*}
L=\left\{ \begin{pmatrix} m_1 & 0\\
0 & m_2\end{pmatrix} \; \Big|\; m_i\in\GL(n_{i},\R),\;
\det m_1\cdot \det m_2 =1\right\}.
\end{equation*}
So $L$ is isomorphic to $\GL(2,\R)$. The Eisenstein series are associated to
Maass cusp forms on $\Gamma(1)\bs\bH$. 
The constant terms of the Eisenstein series  are described in 
\cite[Proposition 10.11.2]{Go}. Let $f$ be a Maass cusp form for $\Gamma(1)$
and let $\Lambda(s,f)$ be the completed $L$-function of $f$ defined by 
(\ref{3.2a}). Then the relevant constant term of the Eisenstein series
 associated to $f$ is given by
\[\frac{\Lambda(s,f)}{\Lambda(1+s,f)}.\]
To proceed one needs a bound similar to (\ref{4.6}). Assume that $\Delta f=
(1/4+r^2)f$. Using the analytic properties of $\Lambda(s,f)$ one can show that
for $T\ge1$ 
\begin{equation}\label{5.7}
\int_{-T}^T\frac{\Lambda^\prime}{\Lambda}(1+it,f)\;dt\ll 
T\log(T+|r|).
\end{equation}
This is the key result that is needed to deal with the contribution of the
Eisenstein series to the trace formula. 

The example demonstrates a general feature of spectral theory on locally 
symmetric spaces.
Harmonic analysis on higher rank spaces requires the knowledge
of the analytic properties of automorphic $L$-functions attached to cusp
forms on lower rank groups. For $\GL(n)$, the corresponding $L$-functions 
are Rankin-Selberg convolutions $L(s,\phi_1\times \phi_2)$ of automorphic
cusp forms on $\GL(n_i)$, $i=1,2$, where $n_1+n_2=n$ (cf. \cite{Bu}, \cite{Go}
for their definition). The analytic properties
of these $L$-functions are well understood so that estimates similar to
(\ref{5.7}) can be established. 
For other groups $\G$ (except for some low dimensional cases) our current 
knowledge of the analytic properties of the corresponding $L$-functions  
 is not sufficient  to prove estimates like (\ref{5.7}). Only partial results
exist \cite{CPS}. This is one of the main obstacles to extend Theorem 
\ref{th5.3} to other groups.

\end{document}